\documentclass{article}

\usepackage{epsf,epsfig,amsfonts,amsgen,amsmath,amstext,amsbsy,amsopn,amsthm,cases,listings,color
}

\usepackage{ebezier,eepic}
\usepackage{color}
\usepackage{multirow}
\usepackage{epstopdf}
\usepackage{graphicx}
\usepackage{pgf,tikz}
\usepackage{mathrsfs}
\usepackage[marginal]{footmisc}
\usepackage{enumitem}
\usepackage[titletoc]{appendix}
\usepackage{booktabs}
\usepackage{url}
\usepackage{subfigure}
\usepackage{mathrsfs}
\usetikzlibrary{arrows}

\usepackage{authblk}

\allowdisplaybreaks[1]

\definecolor{uuuuuu}{rgb}{0.27,0.27,0.27}
\definecolor{sqsqsq}{rgb}{0.1255,0.1255,0.1255}

\newtheorem{definition}{Definition} 
\newtheorem{theorem}[definition]{Theorem}
\newtheorem{lemma}[definition]{Lemma}
\newtheorem{proposition}[definition]{Proposition}

\newenvironment{pf}{\noindent{\bf Proof}}

\usepackage[utf8]{inputenc}

\title{The Linear $q$-Hypergraph Process}
\author[1]{Sayok Chakravarty}
\author[2]{Nicholas Spanier}
\affil[1,2]{Department and Mathematics, Statistics, and Computer Science, University of Illinois Chicago}

\begin{document}

\maketitle
\begin{center}
\textbf{Abstract}
\end{center}
\begin{center}

\begin{flushleft}
     We analyze a random greedy process to construct $q$-uniform linear hypergraphs using the differential equation method. We show for $q=o(\sqrt{\log n})$, that this process yields a hypergraph with $\frac{n(n-1)}{q(q-1)}(1-o(1))$ edges. We also give some bounds for maximal linear hypergraphs.  
\end{flushleft}
\end{center}

\section{Introduction}
\subsection{$F$-free processes and the Differential Equation Method}
The differential equation method for graph processes was popularized by Wormald in 1999 \cite{wormald1999differential}
to analyze random graph processes. The survey \cite{bennett2022gentle} provides an accessible introduction to the differential equation method. A common application of the differential method is the analysis of the $\cal{F}$-free process where $\cal{F}$ is a family of graphs. This is random process which creates a graph $G_i$ on $n$ vertices  by adding edges uniformly at random one at a time so $G_i$ contains no subgraph in the family $\cal{F}$. The case in which $\mathcal{F}$ is a single graph has been studied for graphs including $K_3$ and $K_4$ to give lower bounds on the Ramsey numbers $r(3,t)$ and $r(4,t)$  \cite{bohman2009triangle, bohman2021dynamic,erdos1995size,kim1995ramsey}. 

This paper uses the differential equation method to construct approximate partial Steiner systems. An $(n,q,t)$ \textit{Steiner system} is a family $\mathcal{H} \subset \binom{[n]}{q}$ so that any $t$ subset is contained in exactly one element of $\mathcal{H}$. When each $t$ subset is contained in at most one element in $\mathcal{H}$, the family is called a \textit{partial Steiner system}. It is easy to see that $|\mathcal{H}| \le \frac{\binom{n}{t}}{\binom{q}{t}}$ for a partial Steiner system $\mathcal{H}$ with equality when the family is a Steiner system. When $|\mathcal{H}| = \frac{\binom{n}{t}}{\binom{q}{t}}(1-o(1))$ as $n \rightarrow \infty$ we say $\mathcal{H}$ is an \textit{approximate Steiner system}. The existence of approximate Steiner systems for $q$ constant is proven by R\"{o}dl in \cite{rodl1985packing}. A design with parameters $(n,q,r,\lambda)$ is a collection $S$ of $q$ sets from $[n]$ where each $r$ subset of $[n]$ is contained in exactly $\lambda$ sets in $S$. In the groundbreaking paper \cite{keevash2014existence} Keevash proves that designs exist given the necessary divisibility conditions. Furthermore, in \cite{keevash2018existence} this result was generalized to the setting of subset sums in lattices with coordinates indexed by labelled faces of simplicial complexes. 

Theorem 7.1 of \cite{wormald1999differential} uses the differential equation method to show that a greedy matching of a $k$-uniform hypergraph will use almost all of the vertices given certain degree conditions are satisfied. The problem of finding an $(n,q,t)$ partial Steiner system on $[n]$ can be viewed as finding a matching in a particular $k$-uniform hypergraph. In \cite{wormald1999differential}, Wormald analyzes the greedy packing process to construct a hypergraph matching. He comments that while the proof only works for fixed $k$, one should be able to let $k$ be a function of the number of vertices and get an analogous result. Wormald's result suggests the greedy packing process could construct a $(n,q,2)$ approximate partial Steiner system for $q=o(\sqrt[4]{\log n})$. We explain the connection between Wormald's result and our result in Section 1.2. Bohman, Frieze, and Lubetzky studied a random triangle removal process which constructs a partial Steiner triple system \cite{bohman2015random}. We analyze a process that is equivalent to randomly removing $K_q$ subgraphs. Bennett, Dudek, and Zerbib analyze a similar process they call the online triangle packing process to prove Tuza's conjecture for $G(n,m)$ under certain conditions \cite{bennett2020large}. Our main contribution is that we show that an $\mathcal{F}$-free process constructs an $(n,q,2)$ approximate partial Steiner system for $q=o(\sqrt{\log n})$. 

Bohman and Warnke showed there exists approximate partial Steiner triple systems with high girth by analyzing an $F$-free process \cite{bohman2019large}. Our work uses their approach to analyze the process of $q$-uniform graphs where $q$ may depend on $n$.

\subsection{Notation}
For a sequence of events $\{E_n\}$ we say that this sequence occurs with high probability if 
$$\lim_{n \to \infty} \mathbb{P}(E_n) = 1.$$
We abbreviate with high probability as whp. We will use the notation $f(n) \ll g(n)$ interchangeably with $f(n) = o(g(n))$ to mean $\lim_{n \rightarrow \infty} \frac{f(n)}{g(n)}=0$. For a sequence of random variables $\{X_i\}_{i=1}^n$, let $\Delta X_i = X_{i+1} - X_i$.

\subsection{$q$-Linear Process}
\begin{flushleft}
A hypergraph $\mathcal{H}$ is called \textit{linear} if for any $A,B \in E(\mathcal{H})$ we have $|A \cap B| \le 1$. In other words, any pair of vertices appears in at most one edge. Suppose that $\mathcal{H}$ is $q$-uniform. Then $|E(\mathcal{H})| \le \frac{\binom{n}{2}}{\binom{q}{2}}=\frac{n(n-1)}{q(q-1)}$ as the $\binom{q}{2}$ pairs in each edge are distinct. \\ \vspace{3mm}
Consider the following simple randomized greedy algorithm for constructing a maximal $q$-uniform linear hypergraph where we add one edge at each step. Let $\mathcal{H}_i$ be the hypergraph at step $i$ and let $e_i$ be the edge added at step $i$. \\ \vspace{3mm}
1. Let $\mathcal{H}_0$ be the empty $q$-uniform hypergraph on $[n]$. \\ \vspace{3mm}
2. For $i \ge 1$, in step $i$ pick $e_i$ uniformly at random from the set $$\left\{e \in \binom{[n]}{q} : | e \cap e_j | \le 1 \text{ for } 1 \le j \le i-1 \right\}$$
and form $\mathcal{H}_i$ by adding $e_i$ to $\mathcal{H}_{i-1}$. \\ \vspace{3mm}
We call this the \textit{$q$-linear process}. We analyze the stopping time of this algorithm for various values of $q$ that may depend on $n$. Our first result gives a lower bound on sizes of maximal linear hypergraphs. Note that such lower bounds give a lower bound on the stopping time of the $q$-linear process.
\begin{proposition}
 \label{n^2/q^3}
 Let $\mathcal{H}$ be a $q$-uniform linear hypergraph on $[n]$ that is maximal (i.e. no edge can be added to $\mathcal{H}$ while maintaining linearity). Then $e(\mathcal{H}) \ge \frac{n(n-q+1)}{q(q-1)^2}=\frac{n^2}{q^3}(1-o(1))$.
\end{proposition}

Observe that Proposition \ref{n^2/q^3} shows that any maximal partial $(n,q,2)$ Steiner system asymptotically has at least $\frac{n^2}{q^3}$ edges. In particular the $q$-linear process must continue for at least $\frac{n^2}{q^3}(1-o(1))$ steps.

 In addition to Proposition \ref{n^2/q^3}, notice that there is a trivial lower bound for the size of a maximal $(n,q,2)$ partial Steiner system on the order of $\frac{n}{q}$ by a counting argument. Notice that if $q > \sqrt{n}$, we have that $\frac{n}{q} > \frac{n^2}{q^3}$ and so in this case, the trivial lower is better than the bound from Proposition \ref{n^2/q^3}. Furthermore, when the trivial lower bound is better than the bound from Proposition \ref{n^2/q^3}, Proposition \ref{sqrt(2n)} says exactly how long the process continues asymptotically.

\begin{proposition}
\label{sqrt(2n)}
 Let $\mathcal{H}$ be a $q$-uniform linear hypergraph on $[n]$ with $q \geq \sqrt{2n}$. Then $e(\mathcal{H}) < q$. Further, if $\mathcal{H}$ is maximal then $e(\mathcal{H}) = \Theta(\frac{n}{q})$.
 \end{proposition}

The problem of finding an $(n,q,t)$ partial Steiner system on $[n]$ can be viewed as finding a matching in a $\binom{q}{t}$-uniform hypergraph $H$ where $V(H) = \binom{[n]}{t}$ and for each $S \in \binom{[n]}{q}$ $H$ has an edge which corresponds to all of the $t$-sets in $S$. In \cite{wormald1999differential} Wormald defines the greedy packing process on a hypergraph $H$ as the process which picks an edge from $H$ one at a time uniformly at random and then deletes all the vertices in the chosen edge and continues until there are no edges remaining. We state Wormald's result below:

\begin{theorem}
\label{wormald}
Let $H$ be a $k$-uniform hypergraph with $\nu$ vertices where $k$ is a fixed constant. Assume $\nu < r^C$ for some
constant $C,\delta = o(r^{1/3})$
and $r = o(\nu)$. Also if $d(v)$ is the degree of vertex $v$ in $H$ then assume $|d(v) - r| \leq \delta$. Then for any $\epsilon_0 <\frac{1}{9k(k-1)+3}$ a.a.s. at most $\frac{\nu}{r^{\epsilon_0}}$ vertices remain at the end of the greedy packing process applied to $H$.

\end{theorem}

Based on the connection between partial Steiner systems on $[n]$ and matchings in the hypergraph $H$, the number of unused vertices in the greedy packing process is the number of unused pairs at the end of the $q$-linear process. Note that the correspondences between the greedy packing process and the $q$-linear process is given by $\nu = \binom{n}{2}$ and $k=\binom{q}{2}$. Hence, for the partial Steiner system to have $(1-o(1))\frac{n(n-1)}{q(q-1)}$ edges, Wormald's result suggests that if $k$ were allowed to grow as a function of $\nu$ then we would need that $\frac{\nu}{r^{\epsilon_0}} = o(n^2)$.Then using the fact that $r=o(\nu) = o(n^2)$ and $\epsilon_0 = O(\frac{1}{k^2}) = O(\frac{1}{q^4})$, we would need that $q = o(\sqrt[4]{\log n})$.

Our main result allows us to still get almost all of the edges until $q$ is $o(\sqrt{\log n})$, giving an improvement over the expected result from Wormald 1999 \cite{wormald1999differential}. 

  \begin{theorem}
  \label{main_theorem}
    Let $q = o(\sqrt{\log n})$ and let $\mathcal{H}_n$ be a $q$-uniform hypergraph on $[n]$ obtained from the \textit{$q$-linear process}. Then whp $|E(\mathcal{H}_n)| \geq \frac{n(n-1)}{q(q-1)}(1 - o(1))$.
  \end{theorem}

Note that for $q$ between $\sqrt{\log n}$ and $\sqrt{2n}$ all we know about the $q$-linear process is the lower bound from Proposition \ref{n^2/q^3}. \\ \vspace{3mm}

In Section 2 of this paper, we will prove Proposition \ref{n^2/q^3} and Proposition \ref{sqrt(2n)} and in Section 3 we will prove Theorem \ref{main_theorem}. 
\end{flushleft}

\section{Auxiliary Results}

\begin{flushleft}
We will begin by proving Proposition \ref{n^2/q^3}
\vspace{3mm}
\begin{proof} Let $\mathcal{H}$ be a $q$-uniform linear hypergraph on $[n]$ that is maximal. Consider the graph $G$ on $[n]$ whose edge set is pairs that are not present in any edge of $\mathcal{H}$. Then $G$ is $K_q$-free as a $K_q$ in $G$ corresponds to an edge that can be added to $\mathcal{H}$. Then $\binom{n}{2}= e(\mathcal{H}) \binom{q}{2}+e(G)$. Hence, by Tur\'an's Theorem
\begin{align*}
e(\mathcal{H})
& = \binom{q}{2}^{-1} \left(\binom{n}{2}-e(G) \right) \\
& \ge \binom{q}{2}^{-1} \left(\binom{n}{2}-\left(1-\frac{1}{q-1}\right)\frac{n^2}{2} \right) \\
& = \frac{n(n-q+1)}{q(q-1)^2}.
\end{align*}
\end{proof}
\vspace{3mm}
Next, we prove Proposition \ref{sqrt(2n)}.
\begin{proof} Let $\mathcal{H}$ be a $q$-uniform linear hypergraph on $[n]$. Let $e(\mathcal{H}) = m$ and let $E(\mathcal{H}) = \{e_i : i \in [m]\}$. Now define $\mathcal{H}_i$ as the $q$-uniform hypergraph on $[n]$ with $E(\mathcal{H}_i) = \{e_j : j \in [i]\}$. Define $V_i = [n]\setminus \cup_{j=1}^i e_j$ be the set of vertices not used by any edge in $\mathcal{H}_i$. Now notice that since $\mathcal{H}$ is linear, $|e_{i+1} \cap \cup_{j=1}^i e_j| \leq i$, so $|V_i\setminus V_{i+1}| \geq q-i$. Now let $V_0 = [n]$ and notice that for all $ i \in [m]$
\begin{align*}
|V_i| 
& = n - \sum_{j=1}^i |V_{j-1}\setminus V_{j}| \\
& \le n - \sum_{j=1}^i q - (j-1) \\
& = n - iq + \frac{1}{2}i(i-1).
\end{align*}

Now let $f(x) = n - qx + \frac{1}{2}x(x-1)$ and notice that $f(i) \geq |V_i|$ for all $i \in [m]$. Also, notice that $|V_i|$ is a non-negative integer for all $i \in [m]$ since it is the number of vertices not used in any edge of $\mathcal{H}_i$. Then notice that $f(q) = n - \frac{1}{2}q^2 - \frac{1}{2}q < 0$ since $q \geq \sqrt{2n}$, so if $m \geq q$ this would lead to a contradiction since $|V_q| \leq f(q) < 0$ but $|V_q| \geq 0$. Thus $m < q$.
\\
\vspace{3mm}
Now notice that $\sum_{i=1}^m |V_{i-1} \backslash V_i| \geq \sum_{i=1}^m (q - (i-1)) = mq - \frac{1}{2}m(m-1)$ and further, $\sum_{i=1}^m |V_{i-1} \backslash V_i| \leq n$. Thus we get that $m(q - \frac{1}{2}(m-1)) \leq n$ but since $m-1 < q$ then $m(q - \frac{1}{2}q) \leq n$. Thus $m \leq \frac{2n}{q}$ so $m = O(\frac{n}{q})$.
\\
\vspace{3mm}
Next, assume $\mathcal{H}_i$ is maximal and notice that every new edge uses at most $q$ vertices not used by other edges, and there cannot be $q$ unused vertices because $\mathcal{H}$ is maximal. Thus $m > \frac{n-q}{q} = \Omega(\frac{n}{q})$. Thus $e(\mathcal{H}) = \Theta(\frac{n}{q})$.
\end{proof}
\end{flushleft}

\section{Analysis of the $q$-Linear Process }
\begin{flushleft}
We prove Theorem \ref{main_theorem} using the differential equation method.  
\end{flushleft}
 
\subsection{Trajectories and Definitions}
\begin{flushleft}
To understand $q$-linear process we need to track the codegree of sets $A \subset [n]$. The \emph{codegree} of $A$ at step $i$ is the number of $B \subset [n]$ with $A \cap B=\emptyset$ so that $A \cup B$ can be added to $\mathcal{H}_{i-1}=\left\{e_1,...,e_{i-1}\right\}$. Towards this end, for each $J \subset [n]$ with $|J|=j \in \{0\}\cup [q-1]$ consider the sets
\begin{align*}
& H(i) := \left\{e \in \binom{[n]}{q}: |e \cap e_k| \le 1 \text { for } 1 \le k \le i-1 \right\} \\
& P_j(i) := \left\{J \in \binom{[n]} {j} : J \subset e \text{ for some } e \in H(i)
\right\}\\
& Y_J(i) := \begin{cases}
\left\{K \in \binom{[n] \setminus J}{q-j} : J \cup K  \in H(i)\right\} & J \in P_j(i) \\
Y_J(i-1) & J \not\in P_j(i).
\end{cases} 
\end{align*}
Here $P_j(i)$ represents $j$ sets which can still be a subset of a new edge at step $i$, and $|Y_J(i)|$ represents the codegree of a set in $P_j(i)$ with the convention that if $J \not \in P_j(i)$ then we freeze $Y_J$ at its current value. We are particularly interested in $|Y_{\emptyset}(i)|$ as this gives the number of available edges at step $i$, and we give this set another name $H(i)$ for clarity.
\\
\vspace{3mm}
Next we will define trajectory functions which we expect the random variables $|Y_J(i)|$ to follow. Observe that after $i$ steps the proportion of pairs that are not in any edge is $\frac{\binom{n}{2}-i\binom{q}{2}}{\binom{n}{2}} = 1 - \frac{iq(q-1)}{n(n-1)}$. Our heuristic is that the probability that a pair is not in any edge at step $i$ is  $1 - \frac{iq(q-1)}{n(n-1)}$ and the events that distinct pairs are not in any edges are mutually independent. Now we will define a continuous time variable $t$ which relates to discrete steps by 
$$t(i) = t_i = \frac{i}{n(n-1)}$$
and define the following functions

\begin{align*}
    & p(t) := 1 - q(q-1)t \\
    & y_j(t) := \binom{n-j}{q-j}p^{\binom{q}{2} - \binom{j}{2}} \text{ for all } j \in [q-1] \cup \{0\} \\
    & h(t) := {\binom{n}{q}}p^{\binom{q}{2}}.
\end{align*}

Notice that $p$ is a real-valued function which matches our heuristic for the probability that a pair is not in an edge at step $i$ when $t = \frac{i}{n(n-1)}$. Further, notice that if our heuristic is close to true, then $y_j(t_i)$ gives the approximate size of $Y_J(i)$ when $t = \frac{i}{n(n-1)}$ given that $J \in P_j$. Also note that $h(t) = y_0(t)$.
\\
\vspace{3mm}
Now, we define our targeted stopping time $m_0$ and the errors $\epsilon_j$ we allow on these trajectories as

\begin{align*}
    & f := (\log \log n)^2 \\
    & \beta := \frac{1}{6q^2}\\
    & m_0 := \left \lfloor  \frac{n(n-1)}{q(q-1)}(1-n^{-\beta})\right \rfloor \\
    & \epsilon_j(t) := {\binom{n-j}{q-j}}n^{-1 + 3\beta \binom{q}{2}}q^f p^{-\binom{j}{2} - 2\binom{q}{2}} \\
    & \epsilon_H(t) := \epsilon_0(t).
\end{align*}

Now to prove Theorem \ref{main_theorem}, we prove the following lemma:
\\
\vspace{3mm}
\begin{lemma}
\label{trajectory_lemma}
For all $0 \leq i \leq m_0$, and for all $j \in [q-1] \cup \{0\}$ we have that
\begin{align*}
    & \left||H(i)| - h(t_i) \right| \leq \epsilon_H(t_i) \\
    & ||Y_J(i)| - y_j(t_i)| \leq \epsilon_j(t_i) \textnormal{ for all } J \in P_j(i)
\end{align*}
whp.
\end{lemma}

 We now prove Theorem \ref{main_theorem} assuming Lemma \ref{trajectory_lemma}.

\begin{proof}

Now notice that if $\epsilon_H(t) = o(h(t))$ and $\epsilon_j(t) = o(y_j(t))$, Lemma \ref{trajectory_lemma} will show that whp $H(i) \sim h(t_i)$ and $Y_J(i) \sim y_j(t_i)$, which since $h(t_{m_0}) \gg 1$ will prove Theorem \ref{main_theorem}. 
\\
\vspace{2mm}

To see that $\epsilon_j(t) = o(y_j(t))$ for all $j \in [q-1] \cup \{0\}$ notice that for all $t \in [0, \frac{m_0}{n(n-1)}]$

\begin{align*}
    \frac{\epsilon_j(t)}{y_j(t)}
    & = \frac{\binom{n-j}{q-j}n^{-1+3\beta\binom{q}{2}}q^f p^{-\binom{j}{2}} - 2\binom{q}{2}}{\binom{n-j}{q-j}p^{\binom{q}{2} - \binom{j}{2}}} \\
    & = q^f n^{-1 + 3\beta \binom{q}{2}}p^{-3\binom{q}{2}} \\
    & =  q^f n^{-1 + 3\beta \binom{q}{2}}(1-q(q-1)t)^{-3\binom{q}{2}} \\
    & \le  q^f n^{-1 + 3\beta\binom{q}{2}}\left( 1-q(q-1) \frac{m_0}{n(n-1)}\right)^{-3\binom{q}{2}} \\
    &\leq q^f n^{-1 + 6\beta \binom{q}{2}} = o(1),
\end{align*}

where the last statement in the above follows from our choices of $f$ and $\beta$.
\end{proof}

Notice that since we want $n^{-\beta} = o(1)$ then we need $\beta = \omega\left(\frac{1}{\log n}\right)$, and we also need $\beta = \Theta\left( q^{-2} \right)$. Thus we need that $\binom{q}{2} = o(\log n)$ which holds since we assumed $q = o(\sqrt{\log n})$.
\\
\vspace{2mm}

\end{flushleft}
\subsection{Expected One-Step Change of $Y_J$}
\begin{flushleft}
Let $\Delta Y_{J}(i)=|Y_{J}(i+1)|-|Y_{J}(i)|$ and let $\mathcal{F}_i$ be the natural filtration of the process at step $i$. We refer to $\Delta Y_{J}(i)$ as the one step change of $Y_J(i)$. \\ \vspace{3mm}
\begin{lemma}
\label{One Step change lemma}
The one step change $\mathbb{E}[\Delta Y_J(i) | \mathcal{F}_i]$ is given by
\begin{equation} \label{onestepchange}
  -\frac{1}{|H(i)|} \left(\sum_{K \in Y_J(i)} \left( \sum_{S \subset J, \, T \subset K,\, |S|+|T|\ge 2, \, |T| \ge 1} (-1)^{|S|+|T|} (|S|+|T|-1) |Y_{S \cup T}(i)| \right) \right).  
\end{equation}
\end{lemma}

\begin{proof}
    Observe that $\Delta Y_{J}(i)$ is the number of elements in the codegree of $Y_J(i)$ that are made unavailable by the addition of $e_i$ to $\mathcal{H}_{i-1}$. We show that for fixed $K \in Y_J(i)$ the number of edges $e_i$ that causes $K \notin Y_J(i+1)$ is 
\begin{equation}
\label{One step change inner sum}
\sum_{S \subset J, \, T \subset K,\, |S|+|T|\ge 2, \, |T| \ge 1} (-1)^{|S|+|T|} (|S|+|T|-1) |Y_{S \cup T}(i)|.
\end{equation}
Suppose $e \in H(i)$ with $|e \cap J|=k$ and $|e \cap K|=\ell$ such that $k + \ell \geq 2$. We will use the following Lemma to show that if $k=0,1$ then $e$ is counted once in (\ref{One step change inner sum}) and that if $k \ge 2$ then $e$ is counted 0 times in (\ref{One step change inner sum}).\\ \vspace{3mm}

\begin{lemma}
\label{Combinatorial Identities Lemma}
The following identities hold
$$\sum_{m=2}^\ell \binom{\ell}{m} (-1)^m (m-1)=1$$
$$\sum_{m=2}^\ell \binom{\ell}{m} (-1)^m (m-1)+\sum_{m=1}^\ell \binom{\ell}{m} (-1)^{m+1}m=1$$
$$\sum_{0 \le m_1 \le k, 1 \le m_2 \le \ell, m_1+m_2 \ge 2} \binom{k}{m_1} \binom{\ell}{m_2} (-1)^{m_1+m_2}(m_1+m_2-1)=0.$$
\end{lemma}
We leave the proof of Lemma \ref{Combinatorial Identities Lemma} to the appendix.

Let $k=0$. Then $e$ is counted 
$$\sum_{m=2}^\ell \binom{\ell}{m} (-1)^m (m-1)=1$$
times in (\ref{One step change inner sum}). \\ \vspace{3mm}
Let $k=1$. Then $e$ is counted 
$$\sum_{m=2}^\ell \binom{\ell}{m} (-1)^m (m-1)+\sum_{m=1}^\ell \binom{\ell}{m} (-1)^{m+1}m=1$$
times in (\ref{One step change inner sum}). \\ \vspace{3mm}
Let $k \ge 2$. Then $e$ is counted 
$$\sum_{0 \le m_1 \le k, 1 \le m_2 \le \ell, m_1+m_2 \ge 2} \binom{k}{m_1} \binom{\ell}{m_2} (-1)^{m_1+m_2}(m_1+m_2-1)=0$$
times in (\ref{One step change inner sum}). \\ \vspace{3mm}
\end{proof}
\end{flushleft}

\subsection{Supermartingale and Submartingale Properties}
\begin{flushleft}
Let $\mathcal{G}_i$ be the event that all the bounds
in Lemma \ref{trajectory_lemma} hold for all $j \leq i$. To prove Lemma \ref{trajectory_lemma}, we will define the following random variable where $J \in {\binom{[n]}{j}}$ for $j \in [q-1] \cup \{0\}$ as 
\begin{align*}
    Y_J^{\pm}(i) = \begin{cases}
    |Y_J(i)| - (y_j(t_i) \pm \epsilon_j(t_i)), & \text{if } \mathcal{G}_{i-1} \text{ holds and } J \in P_j(i) \\
    Y_J^{\pm}(i-1), & \text{otherwise}.
    \end{cases}
\end{align*}
\begin{lemma}
    \label{Supermartingale and Submartingale Lemma}
    Let $n \ge n_0$ for some sufficiently large constant $n_0$. For all $J \subseteq [n]$ with $|J| \leq q-1$, $\{Y_J^+(i)\}_{i=0}^{m_0}$ is a supermartingale and $\{Y_J^-(i)\}_{i=0}^{m_0}$ is a submartingale.
\end{lemma}

\begin{proof}
We first note that 
$$\Delta Y_J^{+}(i) = \left(Y_J(i+1)-Y_J(i) \right)-(y_j(t_{i+1})-y_j(t_i))-(\epsilon_j(t_{i+1})-\epsilon_j(t_i)).$$
Since by Taylor's theorem
$$y_j(t_{i+1})-y_j(t_i)=\frac{y_j^{\prime}(t_i)}{n(n-1)}+\frac{1}{2} \frac{y_j^{\prime \prime}(c)}{n^2(n-1)^2}$$
for some $c \in [t_i,t_{i+1} ]$, and similarly 

$$\epsilon_j(t_{i+1})-\epsilon_j(t_i)=\frac{\epsilon_j^{\prime}(t_i)}{n(n-1)}+\frac{1}{2} \frac{\epsilon_j^{\prime \prime}(c)}{n^2(n-1)^2}$$
for some $c \in [t_i,t_{i+1} ]$, we have that $\mathbb{E}[\Delta Y_J^+(i) | \mathcal{F}_i]$ is at most
\begin{equation} \label{changeY_J^+}
\mathbb{E}[\Delta Y_J(i) | \mathcal{F}_i]-\frac{y_j^{\prime}(t_i)}{n(n-1)}-\frac{\epsilon_j^{\prime}(t_i)}{n(n-1)}+\frac{\sup_{s \in [0,\frac{m_0}{n(n-1)}]} |y_j^{\prime \prime}(s)|}{2n^2(n-1)^2}+\frac{\sup_{s \in [0,\frac{m_0}{n(n-1)}]} |\epsilon_j^{\prime \prime}(s)|}{2n^2(n-1)^2}.
\end{equation}

Note that if $\mathcal{G}_{i-1}$ does not hold then $\Delta Y_J^{+}(i)=0$. In the event $\mathcal{G}_{i-1}$, by Lemma \ref{One Step change lemma} we get
\begin{equation}
\label{One step change in G_i event}
\begin{split}
\mathbb{E}[\Delta Y_J(i) | \mathcal{F}_i]
& \le \frac{y_j(t_i)+\epsilon_j(t_i)}{h(t_i)-\epsilon_h(t_i)} \sum_{m=2}^q \left( \binom{q}{m}-\binom{j}{m} \right)  (m-1) ((-1)^{m+1} y_m(t_i)+\epsilon_m(t_i)) \\
& \sim \sum_{m=2}^q \left( \binom{q}{m}-\binom{j}{m} \right)   (m-1) \left( (-1)^{m+1} \frac{y_j(t_i) y_m(t_i)}{h(t_i)}+\frac{y_j(t_i) \epsilon_m(t_i)}{h(t_i)} \right). 
\end{split}
\end{equation}
To show that $\left\{ Y_J^+(i) \right\}$ is a supermartingale, we need to verify that $\mathbb{E}[\Delta Y_J^+(i) | \mathcal{F}_i] \le 0$. We do this by showing that the negative terms in (\ref{changeY_J^+}) are larger than the positive terms. We do this by showing the following lemma. 

\begin{lemma}
\label{Martingale_lemma}
The following hold for all large enough $n$ and $t \in \left[0, m_0/(n(n-1)) \right]$
$$\left(\binom{q}{2} - \binom{j}{2}\right)\frac{y_j(t) y_2(t)}{h(t)} = -\frac{y_j^{\prime}(t)}{n(n-1)}$$
$$\binom{q}{2} \frac{y_j(t) \epsilon_2(t)}{h(t)} \leq \frac{1}{2}\frac{\epsilon_j'(t)}{n(n-1)}$$
$$\binom{q}{m} \frac{y_j(t) y_m(t) }{h(t)} \ll \frac{1}{q} \frac{\epsilon_j'(t)}{n(n-1)} \text{ for } m \ge 3$$
$$\binom{q}{m} \frac{y_j(t) \epsilon_m(t) }{h(t)} \ll \frac{1}{q} \frac{\epsilon_j'(t)}{n(n-1)} \text{ for } m \ge 3$$
$$\frac{\sup_{s \in [0,\frac{m_0}{n(n-1)}]} |y_j^{\prime \prime}(s)|}{2n^2(n-1)^2}+\frac{\sup_{s \in [0,\frac{m_0}{n(n-1)}]} |\epsilon_j^{\prime \prime}(s)|}{2n^2(n-1)^2} \ll \frac{\epsilon_j'(t)}{n(n-1)}.$$
\end{lemma}
We leave the proofs of Lemma \ref{Martingale_lemma} for the appendix. Note that the first equation in Lemma \ref{Martingale_lemma} shows that the $-\left( \binom{q}{2} - \binom{j}{2}\right)\frac{y_j(t_i)y_2(t_i)}{h}$ term from (\ref{One step change in G_i event}) cancels completely in (\ref{changeY_J^+}) with $-\frac{y_j'(t_i)}{n(n-1)}$. The second, third, and fourth equations in Lemma \ref{Martingale_lemma} show that the absolute value of the rest of the terms in (\ref{One step change in G_i event}) is smaller than the absolute value of the $-\frac{\epsilon'_j(t_i)}{n(n-1)}$ term in (\ref{changeY_J^+}). Finally the last term in Lemma \ref{Martingale_lemma} shows the absolute value of the second derivative terms is smaller than the absolute value of the $-\frac{\epsilon'_j(t_i)}{n(n-1)}$ term in (\ref{changeY_J^+}).

For showing $\{Y_J^-(i)\}$ is a submartingale, notice that the main term from $\mathbb{E}[\Delta Y_J(i) | \mathcal{F}_i]$ and the $y_j'(t_i)$ term still cancel, so the $\epsilon_j'(t_i)$ term is still sufficiently larger than all the other terms, but the $\epsilon'(t_i)$ term is positive in $\mathbb{E}[\Delta Y_J^{-}(i) | \mathcal{F}_i]$, so $\mathbb{E}[\Delta Y_J^{-}(i) | \mathcal{F}_i] \geq 0$.

\end{proof}
\end{flushleft}
\subsection{Absolute Bound on One-Step Change}
\begin{flushleft}
\begin{lemma}
\label{Absolute Bound One Step Change Lemma}
    For all $J \subseteq [n]$ with $|J| \leq q-1$ and all $i \in [m_0]$,
    $$|\Delta Y_J^+(i)| \leq (q-1)\binom{n-j}{q-j-1}(1+o(1))$$
    $$|\Delta Y_J^-(i)| \leq (q-1)\binom{n-j}{q-j-1}(1+o(1)).$$

\end{lemma}
\begin{proof}
First, notice that
\begin{equation} \label{one step change Y_J^+}
  |\Delta Y_J^{+}(i)| \leq |\Delta Y_J(i)| + \sup_{t \in [t_i, t_{i+1}]} \frac{|y_j'(t)|}{n(n-1)} + \sup_{t \in [t_i, t_{i+1}]} \frac{|\epsilon_j'(t)|}{n(n-1)}.
\end{equation}

We will bound each of the terms in (\ref{one step change Y_J^+}). We will start by bounding $\Delta |Y_J(i)|$. Notice that since $Y_J$ only changes when an available edge containing $J$ becomes unavailable, then existing edges can only cause the absolute change in $Y_J(i)$ to be smaller since sets that would have been removed from the codegree of $J$ were already not in the codegree of $J$. Thus without loss of generality we may assume $i=0$, and now we will consider three types of edges, $e$, which can be added, edges where $|e \cap J| \geq 2$, edges where $|e \cap J| = 1$, and edges where $|e \cap J| = 0$. First, since we freeze $Y_J$ once $J$ has an intersection with an existing edge of size at least 2, in the case where $|e \cap J| \geq 2$ we get that $\Delta Y_J(0) = 0$.

Next, when $|e \cap J| = 1$ then the number of $f \in Y_J(0) \setminus Y_J(1)$ is the number of edges which contain $J$ and at least one vertex in $e \setminus J$. To upperbound this quantity, we pick one of the $q-1$ vertices in $e \setminus J$ and then just pick the rest of the vertices of $f$ as any $q - j - 1$ vertices in $[n] \setminus J$, which gives an upper bound on $|\Delta Y_J(0)|$ of

\begin{align*}
    |\Delta Y_J(0)| &\leq (q-1)\binom{n-j}{q-j-1}.
\end{align*}

We will now show that $(q-1)\binom{n-j}{q-j-1}$ is the largest term in (\ref{one step change Y_J^+}) using the following standard lower bound 

\begin{equation}
    \label{standard binomial estimate}
    \binom{n-j}{q-j-1} \geq \left(\frac{n-j}{q-j-1}\right)^{q-j-1}
    = \Omega(n^{q-j-1}q^{-q+j+1}).
\end{equation}

When $|e \cap J| = 0$, the only element of $Y_J(0)$ that are not in $Y_J(1)$ are those with at least 2 vertices in $e$. It follows that when $|e \cap J|=0$,

\begin{align*}
    |\Delta Y_J(0)| &\leq \binom{q}{2}\binom{n-j}{q-j-2} \\
    &\leq \frac{1}{2}q^2\left(\frac{(n-j)e}{q-j-2}\right)^{q-j-2}\\
    & = O(n^{q-j-2}e^{q-j-2}q^2) \\
    &= o\left((q-1) \binom{n-j}{q-j-1}\right)
    \end{align*}

where the last line follows from (\ref{standard binomial estimate}). Thus for all $J$ and for all $i$ we have that $|\Delta Y_J(i)| \le (q-1) \binom{n-j}{q-j-1}.$
\\
\vspace{3mm}

Next, we will bound $\sup_{t \in [t_i, t_{i+1}]} \frac{|y_j'(t)|}{n(n-1)}$. Notice that for all $t \in [0, m_0]$

\begin{align*}
    \left| \frac{y_j'(t)}{n(n-1)} \right| &= \left| \frac{\binom{n-j}{q-j}p^{\binom{q}{2} - \binom{j}{2} - 1}(\binom{q}{2} - \binom{j}{2})(-q(q-1))}{n(n-1)} \right| \\
    &= O\left(\left(\frac{n e}{q-j}\right)^{q-j}q^4 n^{-2}\right) \\
    &= O(n^{q-j-2}q^{4}e^{q-j}) \\
    &= o\left((q-1) \binom{n-j}{q-j-1}\right)
\end{align*}
where the last line follows from (\ref{standard binomial estimate}). Finally we will bound $\sup_{t \in [t_i, t_{i+1}]} \frac{|\epsilon_j'(t)|}{n(n-1)}$. Now notice that for all $t \in [0, m_0/n(n-1)]$ we get that

\begin{align*}
    \left|\frac{\epsilon_j'(t)}{n(n-1)} \right| &= \left| \frac{\binom{n-j}{q-j}n^{-1 + 3\beta \binom{q}{2}}q^f p^{-\binom{j}{2} - 2\binom{q}{2} - 1}(-\binom{j}{2} - 2\binom{q}{2})(-q(q-1))}{n(n-1)} \right| \\
    &= O\left(\left(\frac{ne}{q-j}\right)^{q-j}n^{-1 + 3\beta \binom{q}{2}}q^f \left(1-q(q-1)\frac{m_0}{n(n-1)}\right)^{-\binom{j}{2} - 2\binom{q}{2} - 1} n^{-2}q^4\right) \\
    &= O(n^{q-j-3+3\beta \binom{q}{2} + \beta(\binom{j}{2} + 2\binom{q}{2} + 1)}q^{4 + f}e^{q-j}) \\
    &= o\left((q-1)\binom{n-j}{q-j-1}\right)
\end{align*}

where the last line follows from (\ref{standard binomial estimate}). Thus $|\Delta Y_J^+(i)| \leq (q-1){\binom{n-j}{q-j-1}}(1+o(1))$ for all $J$ and $i$. A similar proof bounds $|\Delta Y_J^-(i)|$.
\end{proof}
\end{flushleft}
\subsection{Freedman's Inequality}

To finish the proof of Lemma \ref{trajectory_lemma} we use Freedman's Inequality which we state below \cite{freedman1975tail}. \\ \vspace{3mm}
\begin{theorem}
\label{freedman}
 Let $\left\{S(i) \right\}_{i \ge 0}$ be a supermartingale with respect to the filtration $\mathcal{F}=\left\{ \mathcal{F}_i \right\}_{i \ge 0}$. If $\max_{i \ge 0} |\Delta S(i)| \le C$ and $\sum_{i \ge 0} \mathbb{E} ( |\Delta S(i)| \, | \, \mathcal{F}_i) \le V$, then for any $z >0$
$$\mathbb{P} \left( S(i) \ge S(0)+z \text{ for some } i \ge 0 \right) \le \exp \left\{ -\frac{z^2}{2C(V+z)} \right\}.$$
\end{theorem}
We now prove Lemma \ref{trajectory_lemma}.
\begin{proof}
We first apply Theorem \ref{freedman} to $\{Y_J^+\}$ to prove the upper bound in Lemma \ref{trajectory_lemma}. Observe that if we set $S=Y_J^+$ and $z_j=-Y_J^+(0)=\epsilon_j(0)$ and show that $\frac{z_j^2}{2C_j(V_j+z_j)} \rightarrow \infty$ as $n \rightarrow \infty$ then we will have shown that $\mathbb{P}(Y_J^+(i) < 0 \text{ for all }i )$ goes to $1$ as $n$ goes to infinity. This along with the analogous statement of $Y_J^-$ and a union bound argument will show that the inequalities in Lemma \ref{trajectory_lemma} hold. \\ \vspace{3mm}
We now compute $C_j$ and $V_j$. From the absolute bound on the one step change in $Y_J^+$ we know that 
\begin{align*}
|\Delta Y_J^+(i)| 
& \le \binom{q-1}{1} \binom{n-j}{q-j-1}(1 + o(1)).
\end{align*}
So we can take $C_j=\binom{q-1}{1} \binom{n-j}{q-j-1}(1+o(1))$
Furthermore, Lemma \ref{Martingale_lemma} together with $m_0 \leq \frac{n(n-1)}{q(q-1)}$ implies that 
\begin{align*}
\sum_{i \ge 0} \mathbb{E} ( |\Delta Y_J^+(i)| \, | \, \mathcal{F}_i) 
& \le \frac{n(n-1)}{q(q-1)} O \left( \frac{\sup_{t \in [0, t_{m_0}]}|\epsilon_j^{\prime}(t)|}{n(n-1)} \right) \\
& \le \frac{n(n-1)}{q(q-1)} O \left(\left| \frac{\binom{n-j}{q-j}n^{-1 + 3\beta \binom{q}{2}}q^f n^{\beta \left(\binom{j}{2}+ 2\binom{q}{2} + 1 \right)}(-\binom{j}{2} - 2\binom{q}{2})(-q(q-1))}{n(n-1)} \right| \right) \\
& = O \left( \binom{n-j}{q-j} q^{f+2} n^{-1+6\beta \binom{q}{2}} \right).
\end{align*}
Hence, we can set $V_j=O \left( \binom{n-j}{q-j} q^{f+2} n^{-1+6\beta \binom{q}{2}} \right)$. \\ \vspace{2mm}
Set $z_j=\epsilon_j(0)=\binom{n-j}{q-j}n^{-1 + 3\beta \binom{q}{2}}q^f $. Notice that $z_j \ll V_j$, so to verify that $\frac{z_j^2}{2C_j(V_j+z_j)} \gg 1$, it suffices to check that $\frac{z_j^2}{C_jV_j} \gg 1$. Observe that 
\begin{align*}
\frac{z_j^2}{C_jV_j}
& = \frac{\binom{n-j}{q-j}^2n^{-2 + 6\beta \binom{q}{2}}q^{2f}}{O \left( \binom{n-j}{q-j} q^{f+2} n^{-1+6\beta \binom{q}{2}} \right) \cdot \binom{q-1}{1} \binom{n-j}{q-j-1}} \\
& = \Omega \left( \frac{n-q+1}{q-j}  n^{-1}q^{f-3}\right) \\
& = \Omega \left( \left(1-\frac{q}{n}+\frac{1}{n} \right)q^{f-4} \right) \\
& = \Omega(q^{f-4}).
\end{align*} 
By Theorem \ref{freedman}, for all $J$
$$\mathbb{P}(Y_J^+(i) \geq 0 \text{ for some } i) \leq \exp \left(-\Omega(q^{f-4}) \right).$$
Since $\{Y_J^-\}$ is a submartingale, $\{-Y_J^-\}$ is a supermartingale, so by applying Theorem \ref{freedman} to $\{-Y_J^-\}$ with $S = -Y_J^-$ and $z=Y_J^-(0) = \epsilon_j(0)$, a similar argument shows for all $J$ $$\mathbb{P}(-Y_J^-(i) \geq 0 \text{ for some } i) \leq \exp \left(-\Omega(q^{f-4}) \right).$$

To get the conclusion of Lemma \ref{trajectory_lemma}, we show $\mathbb{P}(\mathcal{G}_{m_0}^c) \rightarrow 0$ as $n \rightarrow \infty$. Observe that 
\begin{align*}
\mathbb{P}(\mathcal{G}_{m_0}^c)
& \le \mathbb{P} \left( \bigcup_{J \subset [n], |J|<q} \left\{ Y_J^{+}(i) \ge 0 \text{ for some } i \ge 0\right\} \cup \left\{Y_J^{-}(i) \le 0 \text{ for some } i \ge 0\right\}\right) \\
& \le 2q \binom{n}{q} e^{-\Omega(q^{f-4})} \\
& \le 2q n^q e^{-\Omega(q^{f-4})}.
\end{align*}
Observe that 
$$\log\left(2qn^qe^{-\Omega(q^{f-4})}\right)=\log 2+\log q+ q \log n-\Omega(q^{f-4}).$$
Since the largest term is $q \log n$ we need to verify that that $q \log n \ll q^{f-4}$. To see this note that $f \gg \frac{\log \log n}{\log q}$. Hence, $2q n^q e^{-\Omega(q^{f-4})}  =o(1)$ and we have $\mathbb{P}(\mathcal{G}_{m_0}) \rightarrow 1$ as $n \rightarrow \infty$ which proves Lemma \ref{trajectory_lemma}. 
\end{proof}
\subsection{Appendix}
\subsubsection{Proof of Lemma \ref{Martingale_lemma}}
\begin{flushleft}
We now prove Lemma \ref{Martingale_lemma}. \\ \vspace{3mm}
\begin{proof}
    Throughout this section, we will use the derivatives $y'_j(t)$ and $\epsilon'_j(t)$ which are given by $$y'_j(t) = -\left(\binom{q}{2} - \binom{j}{2} \right)\binom{n-j}{q-j}p^{\binom{q}{2} - \binom{j}{2} - 1}q(q-1)$$ $$\epsilon'_j(t) = \left(\binom{j}{2} + 2\binom{q}{2} \right)\binom{n-j}{q-j}n^{-1 + 3\beta \binom{q}{2}} q^f p^{-\binom{j}{2} -2\binom{q}{2} -1}q(q-1).$$

    We first show that $\left(\binom{q}{2}-\binom{j}{2} \right)\frac{y_j(t)y_2(t)}{h(t)} =- \frac{y_j'(t)}{n(n-1)}$. Observe that
\begin{align*}
    \left(\binom{q}{2} - \binom{j}{2}\right)\frac{y_j(t)y_2(t)}{h(t)} &= \left(\binom{q}{2}-\binom{j}{2}\right)\frac{\binom{n-j}{q-j}p^{\binom{q}{2}-\binom{j}{2}}\binom{n-2}{q-2}p^{\binom{q}{2}-1}}{\binom{n}{q}p^{\binom{q}{2}}} \\
    &= \left(\binom{q}{2} - \binom{j}{2} \right)\binom{n-j}{q-j} \frac{q(q-1)}{n(n-1)}p^{\binom{q}{2}-\binom{j}{2} - 1} \\
    &= -\frac{y_j'(t)}{n(n-1)}.
\end{align*}

Next, we show that $\binom{q}{2} \frac{y_j(t) \epsilon_2(t)}{h(t)} \leq \frac{1}{2} \frac{\epsilon_j'(t)}{n(n-1)}$. Indeed
\begin{align*}
    \binom{q}{2} \frac{y_j(t) \epsilon_2(t)}{h(t)}\left( \frac{\epsilon_j'(t)}{n(n-1)}\right)^{-1} &= \frac{\frac{\binom{q}{2} \binom{n-j}{q-j}  p^{\binom{q}{2} - \binom{j}{2}} \binom{n-2}{q-2} n^{-1 + 3\beta \binom{q}{2}}q^f p^{-\binom{2}{2} - 2\binom{q}{2}}}{\binom{n}{q} p^{\binom{q}{2}}}}{\frac{\binom{n-j}{q-j} n^{-1 + 3\beta \binom{q}{2}} q^f p^{-\binom{j}{2} - 2\binom{q}{2} -1}(\binom{j}{2} + 2\binom{q}{2})(q(q-1))}{n(n-1)}} \\
    &= \frac{\binom{q}{2}}{\binom{j}{2} + 2\binom{q}{2}} \leq \frac{1}{2}.
\end{align*}

To prove the rest of Lemma \ref{Martingale_lemma}, we will first give a lower bound on $\frac{\epsilon_j'(t)}{n(n-1)}$ and then prove this is asymptotically larger than all the remaining terms. First notice that
\begin{align*}
    \frac{\epsilon_j'(t)}{n(n-1)} &= \frac{\binom{n-j}{q-j}n^{-1 + 3\beta \binom{q}{2}}q^f p^{-\binom{j}{2} - 2\binom{q}{2} - 1}(\binom{j}{2} + 2\binom{q}{2})(q(q-1))}{n(n-1)} \\
    &= \Omega \left(\binom{n-j}{q-j} n^{-3 + 3\beta \binom{q}{2}}q^{f + 4}p^{-\binom{j}{2} - 2\binom{q}{2} - 1}\right).
\end{align*}

Now, we will compute $q \binom{q}{m} \frac{y_j(t) y_m(t)}{h(t)}$ where $3 \leq m \leq q-1$ and show that each of these terms is $o(\binom{n-j}{q-j} n^{-3 + 3\beta \binom{q}{2}}q^{f + 4}p^{-\binom{j}{2} - 2\binom{q}{2} - 1})$. Observe that 

\begin{align*}
    q \binom{q}{m} \frac{y_j(t) y_m(t)}{h(t)} &= q \binom{q}{m} \frac{\binom{n-j}{q-j}p^{\binom{q}{2} - \binom{j}{2}} \binom{n-m}{q-m} p^{\binom{q}{2} - \binom{m}{2}}}{\binom{n}{q} p^{\binom{q}{2}}} \\
    &\leq \binom{n-j}{q-j} q \left(\frac{qe}{m}\right)^m \left(\frac{ne}{q-m}\right)^{q-m} \left(\frac{q}{n}\right)^q p^{-\binom{q}{2}} \\
    &\leq \binom{n-j}{q-j}n^{-m + \beta \binom{q}{2}}q^{2m+1}e^{q}m^{-m} \\
    &= o\left(\binom{n-j}{q-j} n^{-3 + 3\beta \binom{q}{2}}q^{f + 4}p^{-\binom{j}{2} - 2\binom{q}{2} - 1}\right).
\end{align*}

Since $\epsilon_m(t) = o(y_m(t))$ 
$$q \binom{q}{m} \frac{y_j(t) \epsilon_m(t)}{h(t)} = o\left(q \binom{q}{m} \frac{y_j(t) y_m(t)}{h(t)}\right) = o\left(\binom{n-j}{q-j} n^{-3 + 3\beta \binom{q}{2}}q^{f + 4}p^{-\binom{j}{2} - 2\binom{q}{2} - 1}\right).$$

Lastly, we need to verify that $\frac{\sup_{s \in \left[0, m_0/(n(n-1))\right]} |y_j''(s)|}{2n^2(n-1)^2}$ and $\frac{\sup_{s \in \left[0, m_0/(n(n-1))\right]} |\epsilon_j''(s)|}{2n^2(n-1)^2}$ are also both $o(\binom{n-j}{q-j} n^{-3 + 3\beta \binom{q}{2}}q^{f + 4}p^{-\binom{j}{2} - 2\binom{q}{2} - 1})$.

\begin{align*}
    \frac{\sup_{s \in [0, \frac{m_0}{n(n-1)}]} |y_j''(s)|}{2n^2(n-1)^2} &= \frac{\binom{n-j}{q-j}(p(0))^{\binom{q}{2} - \binom{j}{2} - 2}(\binom{q}{2} - \binom{j}{2})(\binom{q}{2} - \binom{j}{2} - 1)(q^2(q-1)^2)}{2n^2(n-1)^2} \\
    &= O\left(\binom{n-j}{q-j}n^{-4}q^{8}\right) \\
    &= o\left(\binom{n-j}{q-j} n^{-3 + 3\beta \binom{q}{2}}q^{f + 4}p^{-\binom{j}{2} - 2\binom{q}{2} - 1}\right).
\end{align*}

Similarly, we compute

\begin{align*}
    \frac{\sup_{s \in [0, \frac{m_0}{n(n-1)}]} |\epsilon_j''(s)|}{2n^2(n-1)^2} &\leq O \left( \frac{\binom{n-j}{q-j}n^{-1 + 3\beta \binom{q}{2}} q^{f + 8}\left(p\left(\frac{m_0}{n(n-1)}\right)\right)^{-\binom{j}{2} - 2\binom{q}{2} - 2}}{n^2(n-1)^2} \right) \\
    &= O\left(\binom{n-j}{q-j} n^{-5 + \beta(6 \binom{q}{2} + 2)}q^{f + 8}\right) \\
    &= o\left(\binom{n-j}{q-j} n^{-3 + 3\beta \binom{q}{2}}q^{f + 4}p^{-\binom{j}{2} - 2\binom{q}{2} - 1}\right).
\end{align*}

This completes the proof of Lemma \ref{Martingale_lemma}.
\end{proof} 

\end{flushleft}
\subsubsection{Proof of Combinatorial Identities}
\begin{flushleft}

We now prove Lemma \ref{Combinatorial Identities Lemma}
\begin{proof}
We first show that $\sum_{m=2}^\ell \binom{\ell}{m} (-1)^m (m-1)=1$. Observe that 
$$\frac{(1+x)^{\ell}-1}{x}=\sum_{m=1}^\ell \binom{\ell}{m} x^{m-1}.$$
This means
$$\frac{d}{dx} \left[\frac{(1+x)^{\ell}-1}{x} \right]=\sum_{m=2}^\ell \binom{\ell}{m}(m-1) x^{m-2}.$$
Letting $x=-1$ yields $\sum_{m=2}^\ell \binom{\ell}{m} (-1)^m (m-1)=1$. \\ \vspace{3mm}
We now show that $\sum_{m=2}^\ell \binom{\ell}{m} (-1)^m (m-1)+\sum_{m=1}^\ell \binom{\ell}{m} (-1)^{m+1}m=1$. Observe that 
$$(1+x)^{\ell}=\sum_{m=0}^\ell \binom{\ell}{m} x^{m}.$$
Differentiating and letting $x=-1$ yields $\sum_{m=1}^\ell \binom{\ell}{m} (-1)^{m+1}m=0$, so we have $\sum_{m=2}^\ell \binom{\ell}{m} (-1)^m (m-1)+\sum_{m=1}^\ell \binom{\ell}{m} (-1)^{m+1}m=1$.   \\ \vspace{3mm}
We now show that 
$$\sum_{0 \le m_1 \le k, 1 \le m_2 \le \ell, m_1+m_2 \ge 2} \binom{k}{m_1} \binom{\ell}{m_2} (-1)^{m_1+m_2}(m_1+m_2-1)=0.$$
We show this by verifying that 
$$\sum_{0 \le m_1 \le k, 0 \le m_2 \le \ell, m_1+m_2 \ge 2} \binom{k}{m_1} \binom{\ell}{m_2} (-1)^{m_1+m_2}(m_1+m_2-1)=1$$ 
and the $m_2=0$ part 
$$\sum_{0 \le m_1 \le k, m_1 \ge 2} \binom{k}{m_1} \binom{\ell}{0} (-1)^{0+m_1}(0+m_1-1)=1.$$
This second equality is the same as the first identity we proved. For the first one, observe that 
$$\frac{(1+x)^{k+\ell}-1}{x}=\frac{(1+x)^{k}(1+x)^{\ell}-1}{x}=\sum_{0 \le m_1 \le k, 0 \le m_2 \le \ell, m_1+m_2 \ge 1} \binom{k}{m_1} \binom{\ell}{m_2} x^{m_1+m_2-1}.$$
This implies 
$$\frac{d}{dx} \left[\frac{(1+x)^{k+\ell}-1}{x} \right]=\sum_{0 \le m_1 \le k, 0 \le m_2 \le \ell, m_1+m_2 \ge 2} \binom{k}{m_1} \binom{\ell}{m_2} (m_1+m_2-1)x^{m_1+m_2-2}.$$
Substituting $x=-1$ yields 
$$\sum_{0 \le m_1 \le k, 0 \le m_2 \le \ell, m_1+m_2 \ge 2} \binom{k}{m_1} \binom{\ell}{m_2} (-1)^{m_1+m_2}(m_1+m_2-1)=1.$$   
\end{proof}
\end{flushleft}

\section{Acknowledgements}
We would like to thank Dhruv Mubayi for the ideas used in Proposition \ref{sqrt(2n)} and his guidance throughout. This work was funded by NSF DMS-1763317.
We would also like to thank the referees for their many helpful comments.
\bibliographystyle{plain}
\bibliography{main1111}
\end{document}